\documentclass[journal]{IEEEtran}
\usepackage{amsfonts}
\usepackage{amsthm,amsmath, amssymb}
\usepackage{graphicx}
\usepackage{verbatim}
\usepackage{bbm}
\usepackage{rotating}
\usepackage{array}

\usepackage{algorithmic}
\usepackage{algorithm}

\usepackage{flushend}

\newcommand{\bo}[1]{\boldsymbol{#1}}

\newtheorem{theorem}{Theorem}

\newtheorem{corollary}{Corollary}
\newtheorem{lemma}{Lemma}

\newcommand{\bmat}{\begin{pmatrix}}
\newcommand{\emat}{\end{pmatrix}}
\usepackage{color}

\DeclareMathOperator*{\argmin}{argmin}

\graphicspath{ {./Figures/} }

\begin{document}
%
\title{Asymptotic and bootstrap tests for the dimension of the non-Gaussian subspace}
%
%
%

\author{Klaus~Nordhausen, Hannu~Oja, David~E.~Tyler and Joni~Virta
\thanks{K. Nordhausen, Hannu Oja and J. Virta are with the Department of Mathematics and Statistics, University of
Turku, Turku, FIN-20014, Finland
(e-mail: klaus.nordhausen@utu.fi)}.%
\thanks{D.E. Tyler is with the Department of Statistics, The State University of New Jersey, Piscataway, US.}}

\maketitle

\begin{abstract}
Dimension reduction is often a preliminary step in the analysis of large data sets. The so-called non-Gaussian component analysis searches for a projection onto the non-Gaussian part of the data, and it is then important to know the correct dimension of  the non-Gaussian signal subspace. In this paper we develop asymptotic as well as bootstrap
tests for the dimension based on the popular fourth order blind identification (FOBI) method.


\end{abstract}

\begin{IEEEkeywords}
Fourth order blind identification (FOBI), independent component analysis, non-Gaussian component analysis.
\end{IEEEkeywords}


%
\IEEEpeerreviewmaketitle

\section{Introduction}\label{sec:introduction}

Throughout the paper we assume  the Non-Gaussian Component Analysis (NGCA) model, that is,  $\bo x_1,...,\bo x_n$ is a random sample from the distribution of
\[
\bo x=\bo A \bo z+\bo b
\]
where $E(\bo z)=\bo 0$ and $Cov(\bo z)=\bo I_p$, $\bo A \in \mathbb{R}^{p\times p}$ is nonsingular, $\bo b\in \mathbb{R}^p$, and   $\bo z=(\bo z_1',\bo z_2')'$,  $\bo z_1\in \mathbb{R}^{q}$ and $\bo z_2\in \mathbb{R}^{p-q}$ are independent random vectors, $\bo z_1$ is non-Gaussian and $\bo z_2$ is Gaussian. For non-Gaussian $\bo z_1$,  there is no $\bo a\in \mathbb{R}^q$ such that $\bo a'\bo z_1$ has a normal distribution. For a Gaussian $\bo z_2$,
$\bo a'\bo z_2$ has a normal distribution for all  $\bo a\in \mathbb{R}^{p-q}$.  The idea then is, based on the observations  $\bo x_1,...,\bo x_n$, to make inference on the unknown $q$, $0\le q\le p$, and estimate the non-Gaussian signal and Gaussian noise subspaces determined by $\bo z_1$ and $\bo z_2$, respectively. In the literature it is often preassumed that the dimension  $q$ is known. In this short note we develop tests for  $q$ that are based on the matrix of fourth moments used in fourth order blind identification (FOBI).

The model can also be written as
\[
\bo x=\bo A_1 \bo z_1+\bo A_2 \bo z_2+ \bo b
\]
where now $\bo A_1\in \mathbb{R}^{p\times q}$ and $\bo A_2\in \mathbb{R}^{p\times (p-q)}$ have ranks $q$ and $p-q$, respectively.
The independent random vectors $\bo A_1 \bo z_1$ and $\bo A_2 \bo z_2$  represent the signal and noise parts of $\bo x$.
Note that $\bo A_1$ and $\bo A_2$ are identifiable only up to postmultiplication
by $q\times q$ and $(p-q)\times (p-q)$ orthogonal matrices, respectively.  If $q\ge p-1$ and the components of $\bo z$ are independent, the model is called independent component model,
$\bo A=(\bo A_1,\bo A_2) $ is then identified up to the signs and permutation of its columns. Inference on $\bo A$ or its inverse, the unmixing matrix  $\bo A^{-1}$, is then known as  independent component analysis (ICA). \cite{Risk:2016} and \cite{Virta:2016} assumed independent components but allowed  $1\le q\le p$; we call this approach  non-Gaussian independent component analysis (NGICA).  In our model there is no restriction on the number of Gaussian components and the non-Gaussian signal components can be dependent
of each other. The Gaussian Mixture Models (GMM) with equal covariance matrices and $q+1$ groups  and multivariate skew-normal distributions ($q=1$), for example, are included in this wider model. Our definition of the NGCA model is as in  \cite{Theis:2011} and originally suggested in \cite{Blanchard05}. For recent contributions and references for NGCA
but always with known $q$, see e.g. \cite{Bean2014} and \cite{SasakiNiuSugiyama2016}.

In the independent component analysis (ICA) the fourth-order blind
identification (FOBI) by \cite{Cardoso89} uses the regular covariance matrix
\[
\bo S_1=E\left[(\bo x -E(\bo x)) (\bo x -E(\bo x))'  \right]
\]
and the scatter matrix based on fourth moments
\[
\bo S_2=E\left[r^2 (\bo x -E(\bo x))  (\bo x -E(\bo x))'\right],
\]
where $r^2=(\bo x -E(\bo x))' \bo S_1^{-1} (\bo x -E(\bo x))$,
and finds an unmixing matrix $\bo W\in \mathbb{R}^{p\times p}$ such that
 \[
\bo W \bo S_1 \bo W'=\bo I_p
\ \ \ \ \mbox{and}\ \ \ \
\bo W \bo S_2 \bo W'=\bo D
\]
for some diagonal matrix  $\bo D=diag(d_1,...,d_p)$. $\bo W$ and $\bo D$ then provide the eigenvectors and  eigenvalues of $\bo S_1^{-1}\bo S_2$ and $d_i=E[(\bo W\bo x)_i^4]+p-1$, $i=1,...,p$.  If  in the independent component model  the fourth moments of $z_1,...,z_p$ are distinct, $\bo W$ is uniquely defined up to signs and permutations of its rows and $\bo W\bo x$ has independent components with distributions of $z_1,...,z_p$ (up to signs). For a Gaussian $z_i$, $d_i=p+2$. The matrix $\bo D$ also lists the eigenvalues of the symmetric matrix
\[
\bo R=\bo S_1^{-1/2}\bo S_2\bo S_1^{-1/2}.
\]
Our test constructions for the dimension of the signal space in the wider NGCA model are based  on the estimated eigenvalues of $\bo R$.
The test statistics were proposed already in \cite{Nordhausen:2016} for the  NGICA model but without a careful analysis of their limiting distributions.

The plan in this paper is as follows. In Section~\ref{sec:statistic} our test statistic for testing whether the dimension of the signal space is $k$ is introduced. The asymptotic and bootstrap test versions are provided in Sections~\ref{sec:asy} and ~\ref{sec:bootstrap}, respectively. Also the estimate based on the asymptotic test is discussed. The two strategies are compared in simulations in Section~\ref{sec:simulation} and the paper ends with a discussion on alternative tests.  The proofs are provided in the Appendix.

Throughout the paper we use the following notation. The  first and second moments of the eigenvalues of positive definite and symmetric $\bo R\in\mathbb{R}^{p\times p}$ are denoted by
\[
m_1(\bo R) := tr(\bo R)/p \ \ \mbox{and}\ \ \ \ m_2(\bo R) :=  tr(\bo R^2)/p,
\]
and the variance  of the  eigenvalues is
$
s^2(\bo R):= m_2(\bo R)-m_1^2(\bo R)
$.
{We write  $\mathcal{O}^{p \times k}$ for the set of $p\times k$ matrices with orthonormal columns, $k \le p$.}

\section{Test statistic for the dimension}\label{sec:statistic}

Recall that $\bo x_1,...,\bo x_n$ is a random sample from a distribution of
$
\bo x=\bo A\bo z +\bo b
$
where  $\bo A\in\mathbb{R}^{p\times p}$  is non-singular, $\bo b\in\mathbb{R}^p$,  $\mathbb{E}(\bo z)=\bo 0$ and $Cov(\bo z)=\bo I_p$.
Further $\bo z=(\bo z_1',\bo z_2')'$ where $\bo z_1$ and $\bo z_2$ are independent, $\bo z_1\in \mathbb{R}^q$ is non-Gaussian and $\bo z_2\in \mathbb{R}^{p-q}$ is Gaussian.
We also need to assume that the fourth moments exist and $E[(\bo u'\bo z_1)^4]\ne 3$ for all $\bo u'\bo u=1$, $\bo u\in  \mathbb{R}^q$. This is a bit stronger assumption than non-Gaussianity of $\bo z_1$. With unknown $q$, we then wish to test the null hypothesis
\[
H_{0,k}: \ \ \mbox{exactly $p-k$ eigenvalues of $\bo R$ equal $p+2$  }
\]
stating that the dimension of the signal space is $k$.

Natural estimates of $\bo S_1$, $\bo S_2$ and $\bo R$ are
$$ \widehat{\bo S}_1=\frac 1n \sum_{i=1}^n (\bo x_i-\bar{\bo x}) (\bo x_i-\bar{\bo x})',
$$
$$
\widehat{\bo S}_2=\frac 1n \sum_{i=1}^n    (\bo x_i-\bar{\bo x})(\bo x_i-\bar{\bo x})' \widehat{\bo S}_1^{-1} (\bo x_i-\bar{\bo x})  (\bo x_i-\bar{\bo x})'
 $$
and
$
\widehat{\bo R}=\widehat{\bo S}_1^{-1/2} \widehat{\bo S}_2\widehat{\bo S}_1^{-1/2}
$, respectively.
To test the null hypothesis $H_{0,k}$, we use the test statistic
\[
T_{k}= \min_{\bo U\in \mathcal{O}^{p\times (p-k)}}  m_2\left(\bo U' (\widehat{\bo R}-(p+2)\bo I_p) \bo U\right).
\]
If
\[
\widehat{\bo U}_k=\argmin_{\bo U\in \mathcal{O}^{p\times (p-k)}}  m_2\left(\bo U' (\widehat{\bo R}-(p+2)\bo I_p) \bo U\right)
\]
then $T_k =m_2\left(\widehat{\bo U}_k' (\widehat{\bo R}-(p+2)\bo I_p) \widehat{\bo U}_k\right)$.

Note that \cite{Kankainen07} used  $m_2\left(\widehat{\bo R}-(p+2)\bo I_p \right)$ to test for (full) multivariate normality which is a special case here.
Further note that the estimated projections  (with respect to Mahalanobis inner product) to the noise and signal subspaces are given by
$\widehat{\bo Q}_k= \widehat{\bo S}_1^{1/2}\widehat{\bo U}_k\widehat{\bo U}_k' \widehat{\bo S}_1^{-1/2}$ and
$\widehat{\bo P}_k=\bo I_p-\widehat{\bo Q}_k$, respectively.


\section{Asymptotic test for dimension}\label{sec:asy}

As the test statistics $T_k$ are invariant under affine transformations $\bo A\bo x_i+\bo b$, $i=1,...,n$,
we can  without loss of generality assume that $\bo A=\bo I_p$ and  $\bo b=\bo 0$.
In the following, we consider the limiting behavior of $n(p-k)T_k$, $k=0,...,p-1$, under true  $H_{0,q}$. \\

For $k=0,1,...,p-1$, write $\bo V_k=(\bo 0,\bo I_{p-k})'\in \mathcal{O}^{p\times (p-k)}$ and
\[
T_k^*=
m_2\left(\bo V_k' (\widehat{\bo R}-(p+2)\bo I_p) \bo V_k\right)
\]
Then $T_k\le T_k^*$, $k=0,...,p-1$, and, for the true value $q$,  $nT_q=nT_q^*+o_P(1)$, see
the Appendix. 
We then have the following.

\begin{theorem}\label{Th:main}
Under the previously stated assumptions and under $H_{0q}$,\\
(i) for $k<q$, $T_k\to_P c$ for some $c>0$, \\
(ii) for $k=q$,
$n (p-k) T_k\to_d C_k$,
and \\
(iii) for $k>q$,
$  n (p-k) T_k \le n(p-k) T_k^* \to_d C_k$
\\
where
$$ C_k\sim {2\sigma_1} \chi^2_{(p-k-1)(p-k+2)/2}+ \left({{2\sigma_1}+ {\sigma_2} (p-k)} \right) \chi^2_1$$
with independent chi squared variables $\chi^2_{(p-k-1)(p-k+2)/2}$ and  $\chi^2_1$
and $\sigma_1=Var\left(\|\bo z  \|^2\right)+8$ and $\sigma_2=4$.
\end{theorem}

In the regular testing procedure, the null hypothesis $H_{0,k}$ is the rejected if
\[
n(p-k)T_k \ge c_{k,\alpha}
\]
where the critical point $c_{k,\alpha}$ is determined by $\mathbb{P}(C_k\ge c_{k,\alpha})=\alpha$. Note that the test, although constructed for $H_{0,k}$ can  actually be viewed as a
 consistent  test also for
 $H^*_{0,k}$: At least $p-k$ eigenvalues of $\bo R$ equal $p+2$.
To find $c_{k,\alpha}$ in practice, $\sigma_1$ must be replaced by its consistent estimate $\hat\sigma_1$. Then, with increasing $n$,  (i) for $k<q$, the power of the test
for $H_{0,k}$ goes to one, (ii) for $k=q$, the size of the test for $H_{0,q}$  goes to prespecified $\alpha$, and (iii) for $k>q$ the rejection probability for $H_{0,k}$  tends to be smaller than $\alpha$.

We next discuss the estimation of $\sigma_1$. Write  $\hat {\bo z}_i=\widehat{\bo W}\bo x_i$, $i=1,...,n$.
Then, even without knowing the true value of $q$, the parameter $\sigma_1$ can be consistently estimated by
{$
\hat\sigma_1=\frac 1{n} \sum_{i=1}^n \|\hat {\bo z}_i \|^4-p^2+8
$.}
Note that in the independent component model, we simply have $\sigma_1=\sum_{k=1}^p E(z_k^4)-p+8$
which yields  as estimate $\hat\sigma_1=\frac 1n \sum_{i=1}^n\sum_{k=1}^p (\hat z_i)_k^4-p+8$.
With known $q$, the estimates can be simplified further.

A consistent estimate $\hat q$ of the unknown dimension  $q$ can   be based on the test statistics $T_k$, $k=0,...,p-1$ as follows.
\begin{corollary}
For all $k=0,...,p-1$,  let
  $(c_{k,n})$  be a sequence such that
$c_{k,n}\to \infty$ and $\frac{c_{k,n}}n\to 0$ as $n\to\infty$.
Then
\[
\mathbb{P}(n (p-k) T_k\ge c_{k,n})\ \rightarrow\
\left\{
       \begin{array}{ll}
         1, & \hbox{if $k<q$ ;} \\
         0, & \hbox{if $k\ge q$.}
       \end{array}
     \right.
\]
and
\[
\hat q=\min \{k\ :\ n (p-k)T_k < c_{k,n}\}\to_P q.
\]
\end{corollary}

\section{Bootstrap test for dimension}\label{sec:bootstrap}

 Theorem~\ref{Th:main} shows  that the limiting distribution of $n(p-q)T_q$ (with estimated noise subspace)  and $n(p-q)T_q^*$ (with known noise subspace) are
  the same. This means that, if one uses the asymptotic test in the small sample case, the variation coming from the estimation of the subspace is ignored.
 We therefore propose that the small sample null distribution of $n(p-k)T_k$ should be  estimated by resampling  from a distribution
 for which  the null hypothesis $H_{0,k}$ is true and  which is as similar as possible to the empirical distribution of $\bo x_1,...,\bo x_n$.
 For this type of bootstrap sampling, we also need estimated projections
$\widehat{\bo Q}_k= \widehat{\bo S}_1^{1/2}\widehat{\bo U}_k\widehat{\bo U}_k' \widehat{\bo S}_1^{-1/2}$ and
$\widehat{\bo P}_k=\bo I_p-\widehat{\bo Q}_k$.
\cite{Nordhausen:2016} suggested  the following procedure for the NGCA model.\\

{\bf Generating a bootstrap sample for  $H_{0,k}$:}\ \ {\it
\begin{enumerate}
\item Starting with centered $\bo X\in \mathbb{R}^{n\times p}$, compute $\widehat{\bo S}_1$, $\widehat{\bo S}_2$, $\widehat{\bo R}$,
$\widehat{\bo U}_k$, $\widehat{\bo Q}_k$  and $\widehat{\bo P}_k$.
 \item Take a bootstrap sample   $\widetilde{\bo X}=(\tilde{\bo x}_1, \ldots, \tilde{\bo x}_n)'$ of size $n$ from  $\bo X$.
\item For the $(p-k)$-dimensional noise space to be gaussian, transform
\[ \bo x_i^*= \widehat{\bo P}_k\tilde{\bo x}_i +\widehat{\bo S}_1^{1/2} \widehat{\bo U}_k \bo o_i ,\ \ i=1,...,n,  \]
and  $\bo o_1, \ldots, \bo o_n$ are iid  from $N_{p-k}(\bo 0,\bo I_{p-k})$.
 \item $\bo X^* =(\bo x_1^*,...,\bo x_n^*)'$.
\end{enumerate}}
In \cite{Nordhausen:2016}  bootstrap sampling for the NGICA model was suggested as well.


Let $T=T(\bo X)$ be a test statistic for $H_{0,k}$ such as  $T_k$. If $\bo X_1^*,...,\bo X_M^*$ are independent bootstrap samples as described above and
$T_i^*=T(\bo X_i^*)$ then the bootstrap $p$-value is given by
\[
\hat p= \frac {\# (T_i^*\ge T)+1} {M+1}.
\]

\section{A simulation study}\label{sec:simulation}

To compare the asymptotic and bootstrap tests
we generate data sets obeying the following three NGCA models. As the test are affine invariant, it is not a restriction to use $\bo A=\bo I_p$ and $\bo b=\bo 0$ in simulations.

\begin{figure}[htb]
\begin{centering}
\includegraphics[width=\columnwidth]{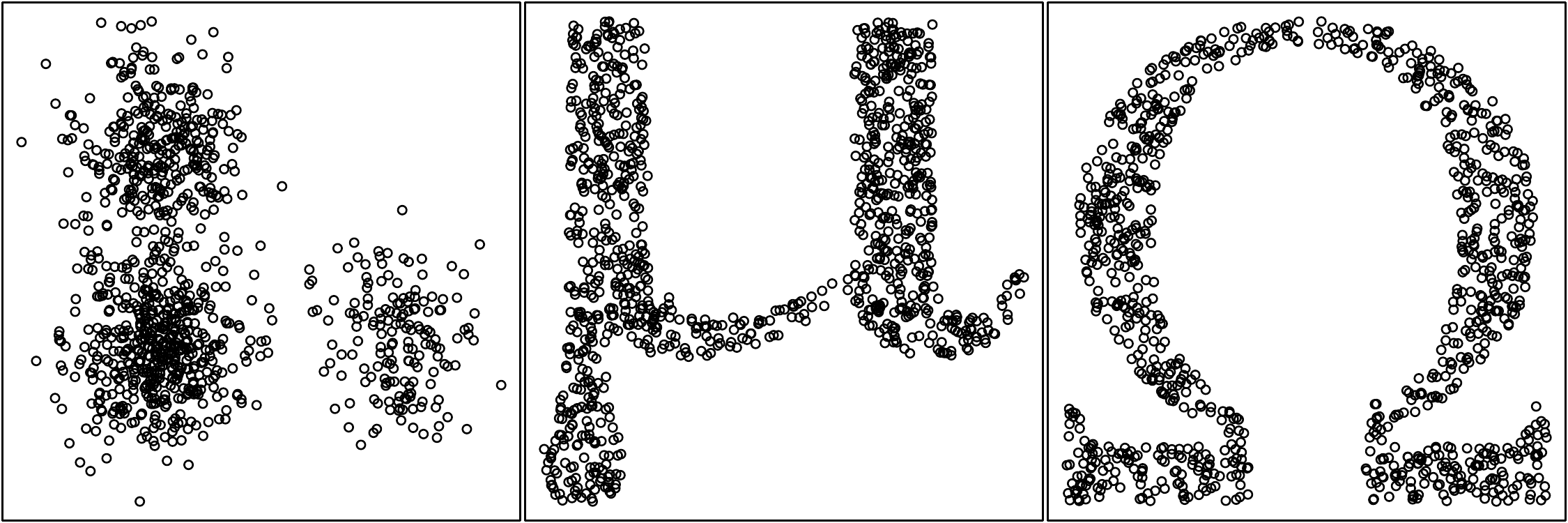}
\caption{Scatterplots for the bivariate non-Gaussian component in model M1 (left) and for two independent bivariate components in model M2.
}
\label{fig1}
\end{centering}
\end{figure}

\begin{description}
  \item[M1:] A GMM model that is a mixture of $N_6(6 \bo e_1, \bo I_6)$, $N_6(4 \bo e_2, \bo I_6)$ and $N_6(\bo 0, \bo I_6)$ with proportions  0.1, 0.4 and 0.5 respectively. Then $p=6$ and $q=2$.
  \item[M2:] An NGCA model with two independent bivariate nongaussian components representing the Greek letters $\bo \mu$ and $\bo \Omega$, see Fig.\ref{fig1}, and independent noise $N_3(\bo 0,\bo I_3)$. Therefore $p=7$ and $q=4$.
  \item[M3:] An NGICA model with the non-Gaussian independent components: exponential, $\chi_1^2$ and uniform and three Gaussian components $N(0,1)$. Hence $p=6$ and $q=3$.
\end{description}

The non-Gaussian components of M1 and M2 are visualized in Figure~\ref{fig1} with samples of size $n=1000$. The complete simulation was performed in R 3.3.2 \cite{R} using the R package ICtest \cite{ICtest} which provides implementations for all methods discussed here. In the comparisons the rejection rates for selected null hypotheses $H_{0,k}$  are reported for asymptotic and bootstrap tests with the test size $\alpha=0.05$.  All results are based on 1000 repetitions.

From the simulations we can conclude that for small sample sizes bootstrapping keeps better the target size 0.05 under the true null hypothesis.  The sample size needed for decent power naturally depends strongly on the underlying model; M1 seems to be  the most difficult case   and requires at least 5000 observations. In general, the results are as suggested by the theory.


{

\begin{table}[ht]
\centering
\small
\caption{Model M1: Rejection rates in 1000 repetitions for   asymptotic and bootstrap tests of  $H_{0,1}$, $H_{0,2}$ (true)  and $H_{0,3}$  with $\alpha=0.05$.}
\label{tab1}
\begin{tabular}{@{\extracolsep{4pt}}lcccccc}
  \hline
   & \multicolumn{2}{c}{$H_{0,1}$} & \multicolumn{2}{c} {$H_{0,2}$} & \multicolumn{2}{c}{$H_{0,3}$}   \\
    \cline{2-3} \cline{4-5} \cline{6-7}
$n$ & Asy & Boot & Asy & Boot & Asy & Boot \\
  \hline
500 &  0.040 & 0.128 & 0.003 & 0.031 & 0.001 & 0.012 \\
  1000 & 0.149 & 0.252 & 0.004 & 0.058 & 0.000 & 0.016 \\
  2000 &  0.403 & 0.466 & 0.027 & 0.068 & 0.000 & 0.024 \\
  5000 &  0.950 & 0.950 &  0.040 & 0.059  & 0.001 & 0.011 \\
  10000 &  0.999 & 0.999  & 0.053 & 0.058 & 0.002 & 0.010 \\
   \hline
\end{tabular}
\end{table}
}

{
\begin{table}[ht]
\centering
\small
\caption{Model M2: Rejection rates in 1000 repetitions for   asymptotic and bootstrap tests of  $H_{0,3}$, $H_{0,4}$ (true)  and $H_{0,5}$  with $\alpha=0.05$.}
\label{tab1}
\begin{tabular}{@{\extracolsep{4pt}}lcccccc}
  \hline
   & \multicolumn{2}{c}{$H_{0,3}$} & \multicolumn{2}{c} {$H_{0,4}$} & \multicolumn{2}{c}{$H_{0,5}$}   \\
    \cline{2-3} \cline{4-5} \cline{6-7}
$n$ & Asy & Boot & Asy & Boot & Asy & Boot \\
  \hline
500 & 0.577 & 0.328 & 0.050 & 0.065 & 0.000 & 0.008 \\
  1000 & 0.997 & 0.964  & 0.073 & 0.056 & 0.002 & 0.008 \\
  2000 &  1.000 & 1.000 & 0.057 & 0.049 & 0.002 & 0.016 \\
  5000 &  1.000 & 1.000 & 0.068 & 0.064 & 0.004 & 0.018 \\
  10000 &  1.000 & 1.000 & 0.050 & 0.046 & 0.050 & 0.046 \\
   \hline
\end{tabular}
\end{table}
}


{
\begin{table}[ht]
\centering
\small
\caption{Model M3: Rejection rates in 1000 repetitions for   asymptotic and bootstrap tests of  $H_{0,2}$, $H_{0,3}$ (true)  and $H_{0,4}$  with $\alpha=0.05$.}
\label{tab1}
\begin{tabular}{@{\extracolsep{4pt}}lcccccc}
  \hline
   & \multicolumn{2}{c}{$H_{0,2}$} & \multicolumn{2}{c} {$H_{0,3}$} & \multicolumn{2}{c}{$H_{0,4}$}   \\
    \cline{2-3} \cline{4-5} \cline{6-7}
$n$ & Asy & Boot & Asy & Boot & Asy & Boot \\
  \hline
500 & 0.715 & 0.796 & 0.024 & 0.051 & 0.001 & 0.015 \\
  1000 &  0.993 & 0.995 & 0.020 & 0.034 & 0.000 & 0.013 \\
  2000 & 1.000 & 1.000 & 0.036 & 0.044 & 0.006 & 0.016 \\
  5000 &  1.000 & 1.000 & 0.042 & 0.041 & 0.003 & 0.011 \\
  10000 &  1.000 & 1.000 & 0.043 & 0.042 & 0.005 & 0.012 \\
   \hline
\end{tabular}
\end{table}

}


\section{Final remarks}\label{sec:final}

In most applications with high-dimensional data only the non-Gaussian variation is informative and the Gaussian part simply presents noise.
Wide literature on NGCA (and NGICA)  provides tools to estimate the non-Gaussian subspace but so far always with known dimension $q$.
In this paper we suggest efficient asymptotic and bootstrap tests for the dimension that are simply based on the eigenvalues of the easily computable FOBI matrix.
Natural estimates of $q$ are also found by successive testing for hypotheses $H_{0,k}$, $k=0,...,p-1$.  Note that the stochastic variability of the $p-k$ eigenvectors depends strongly on how close together the corresponding eigenvalues are. In a similar context \cite{YeWeiss:2003,LuoLi:2016} then suggested a kind of dual estimates of $q$ that are based on  the bootstrap variation of eigenvector estimates.

We end the discussion with some remarks on alternative and competing test statistics. The test statistic $T_k$ can also be written  as a sum  $T_k = T_{k,1}+ T_{k,2}$ where $T_{k,1}= s^2 (\widehat{\bo U}_k' \widehat{\bo R} \widehat{\bo U}_k)$
and $T_{k,2}= [m_1(\widehat{\bo U}_k' \widehat{\bo R} \widehat{\bo U}_k)-(p+2)]^2$.
The first part  $T_{k,1}$ provides a test statistic for the equality of $p-k$  eigenvalues closest to $p+2$   and the second part $T_{k,2}$ measures the deviation
of the average of those eigenvalues from $p+2$ (Gaussian case). Besides $n(p-k)T_k$, one can  also  use
\[
\frac{n(p-k)T_{k,1}}{2\hat\sigma_1} \ \ \mbox{or} \ \
 \frac{n (p-k)T_{k,2}} {2 \hat\sigma_1+4(p-k)}
\]
or their sum  as a test statistic. See the Appendix. Under true  $H_{0,k}$, these statistics have limiting chi square distributions with  $(p-k-1)(p-k+2)/2$, 1,  and $(p-k-1)(p-k+2)/2+1$  degrees of freedom. The first two statistics use less information and are therefore in most cases less powerful than their sum, and the behavior of their sum is very similar to that of $T_k$ as also seen in our simulations (but not reported here).

Note also that FOBI is just a simple special case of the so called two-scatter method \cite{Oja:2006,TCDO09}; $\bo S_1$ and $\bo S_2$ are then replaced by any two scatter matrices
specific to the problem at hand.  Deriving  asymptotic tests  or using bootstrap testing strategy with
\[
T_{k}= \min_{\bo U\in \mathcal{O}^{p\times (p-k)}}  s^2\left(\bo U' \widehat{\bo R} \bo U\right).
\]
for any choices of $\bo S_1$ and $\bo S_2$ is also possible. The properties of these tests with corresponding estimates is a part of our future work.

\section{Appendix: Proof of Theorem~\ref{Th:main} }\label{sec:appendix}

For the limiting distributions of the scatter matrices $\widehat{\bo S}_1$ and $\widehat{\bo S}_2$, we need to assume that the fourth moments of  $\bo z_1$ exist.
Naturally $\bo z_2$ has moments of any order.
Let $\bo x_1,...,\bo x_n$ be a random sample from the distribution of
\[
\bo x=\bo A \bo z+\bo b
\]
where $E(\bo z)=\bo 0$ and $Cov(\bo z)=\bo I_p$,   $\bo z=(\bo z_1',\bo z_2)'$,  $\bo z_1\in \mathbb{R}^{q}$ and $\bo z_2\in \mathbb{R}^{p-q}$ are independent, $\bo z_1$ is non-Gaussian and $\bo z_2$ is Gaussian. Due to affine invariance of the test statistic, it is not a restriction to assume in the following that  $\bo A=\bo I_p$ and $\bo b=\bo 0$.

We write
$ \widehat{\bo S}_1=\frac 1n \sum_{i=1}^n (\bo x_i-\bar{\bo x}) (\bo x_i-\bar{\bo x})'$
and
$
 \widehat{\bo S}_2=\frac 1n \sum_{i=1}^n    (\bo x_i-\bar{\bo x})(\bo x_i-\bar{\bo x})' \widehat{\bo S}^{-1} (\bo x_i-\bar{\bo x})  (\bo x_i-\bar{\bo x})'
$
and, for known $\bo A=\bo I_p$ and $\bo b=\bo 0$, we have
$$ \widetilde{\bo S}_1=\frac 1n \sum_{i=1}^n \bo x_i \bo x_i'
\ \ \mbox{and} \ \
 \widetilde{\bo S}_2=\frac 1n \sum_{i=1}^n   \bo x_i \bo x_i' \bo x_i \bo x_i'
$$
Let
\[
\widehat{\bo R}=\widehat{\bo S}_1^{-1/2} \widehat{\bo S}_2\widehat{\bo S}_1^{-1/2}
\]
be partitioned as
\[
\widehat {\bo R}=\left(
   \begin{array}{cc}
    \widehat{\bo  R}_{11}  & \widehat{\bo  R}_{12} \\
     \widehat{\bo  R}_{21} & \widehat{\bo  R}_{22} \\
   \end{array}
 \right).
\]
{
Theorem~\ref{Th:main} is then
implied by the following four Lemmas, starting with a linearization result for $ \widehat{\bo  R}_{22}\in \mathbb{R}^{(p-q)\times (p-q)}$.
}
\begin{lemma}\label{lemma:R22} Under the stated assumptions,
\begin{eqnarray*}
\sqrt{n}(\widehat{\bo R}-(p+2)\bo I_{p})_{22} &=& \sqrt{n}(\widetilde{\bo R}-(p+2)\bo I_{p})_{22}+o_P(1)
\end{eqnarray*}
where
\begin{eqnarray*}
\widetilde{\bo R}  &:=& \widetilde{\bo S}_2-(p+4)(\widetilde{\bo S}_1-\bo I_p) - \ tr(\widetilde{\bo S}_1-\bo I_p) \bo I_{p}.
\end{eqnarray*}
\end{lemma}

Write $\bo V_k=(\bo 0,\bo I_{p-k})'\in \mathcal{O}^{p\times (p-k)}$.
As, for all $k\ge p$,
\[
\bo U'\bo V_k'(\widetilde{\bo R}-(p+2)\bo I_{p})\bo V_k \bo U \sim \bo V_k'(\widetilde{\bo R}-(p+2)\bo I_{p})\bo V_k\]
for all $\bo U\in \mathcal{O}^{(p-k)\times (p-k)}$ and $\widetilde{\bo R}$ is the average of iid matrices, we  further obtain the following.

\begin{lemma}
Under  the stated assumptions and $k\ge q$,
\[ \sqrt{n}(\bo V_k'(\widehat{\bo R}-(p+2)\bo I_{p})\bo V_k)\to_d \bo N_k\]
where $vec(\bo N_k)$ has a $(p-k)\times (p-k)$-variate normal distribution with zero mean vector and covariance matrix
\[
 \sigma_1 \left(\bo I_{(p-k)^2}+\bo K_{p-k,p-k}\right)+\sigma_2 vec(\bo I_{p-k})vec(\bo I_{p-k})',
\]
 $\bo K_{p,p}=\sum_{i=1}^p\sum_{j=1}^p(\bo e_i\bo e_j')\otimes
(\bo e_j\bo e_i')$ is the commutation matrix, and
\[
\sigma_1=Var\left(\|\bo z \|^2\right)+8 \ \ \mbox{and}\ \ \sigma_2=4.
\]
\end{lemma}

For $k\ge q$, write
\begin{eqnarray*}
T_k^* &=& T_{k,1}^*+ T_{k,2}^*\\
      &:=& s^2 \left({\bo V}_k' \widehat{\bo R} {\bo V}_k\right)+{\left[m_1\left({\bo V}_k' \widehat{\bo R} {\bo V}_k\right)-(p+2)\right]^2}.
\end{eqnarray*}
As in \cite{Nordhausen:2016}, we can show the following
\begin{lemma}\label{Lemma:stars}
Under the stated assumptions  and {$k\ge q$}, the random variables
$nT_{k,1}^*$ and $nT_{k,2}^*$ are asymptotically independent and
\[
\frac{n(p-k)T_{k,1}^*}{2\sigma_1}\to_d \chi^2_{(p-k-1)(p-k+2)/2}
\] and \[
 \frac{n(p-k) T_{k,2}^*} {2 \sigma_1+(p-k)\sigma_2} \to _d \chi^2_1
\]
where
{$\sigma_1=Var\left(\|\bo z \|^2\right)+8$} and $\sigma_2=4$.
\end{lemma}

\begin{lemma}\label{L1:ICA}
Under the stated assumptions, $T_{k}\le T_k^*$, for all $k\ge q$, and
\begin{eqnarray*}
 n  T_q &=& n T_q^* +o_P(1).
\end{eqnarray*}
\end{lemma}

The first part in the Lemma~\ref{L1:ICA} is trivial, the second part follows from Lemma 3.1 in \cite{EatonTyler91}.

\section*{Acknowledgements}
 The research of K. Nordhausen, H. Oja and J. Virta was partially supported by the Academy of
Finland (grant 268703). D.E. Tyler's research was partially supported by the National Science Foundation
Grant No. DMS-1407751. Any opinions, findings and conclusions or recommendations expressed in this material
are those of the author(s) and do not necessarily reflect those of the National Science Foundation.

\clearpage


\bibliographystyle{IEEEtran}
\bibliography{NGCA_references}

\begin{thebibliography}{10}
\providecommand{\url}[1]{#1}
\csname url@samestyle\endcsname
\providecommand{\newblock}{\relax}
\providecommand{\bibinfo}[2]{#2}
\providecommand{\BIBentrySTDinterwordspacing}{\spaceskip=0pt\relax}
\providecommand{\BIBentryALTinterwordstretchfactor}{4}
\providecommand{\BIBentryALTinterwordspacing}{\spaceskip=\fontdimen2\font plus
\BIBentryALTinterwordstretchfactor\fontdimen3\font minus
  \fontdimen4\font\relax}
\providecommand{\BIBforeignlanguage}[2]{{%
\expandafter\ifx\csname l@#1\endcsname\relax
\typeout{** WARNING: IEEEtran.bst: No hyphenation pattern has been}%
\typeout{** loaded for the language `#1'. Using the pattern for}%
\typeout{** the default language instead.}%
\else
\language=\csname l@#1\endcsname
\fi
#2}}
\providecommand{\BIBdecl}{\relax}
\BIBdecl

\bibitem{Risk:2016}
B.~B. Risk, D.~S. Matteson, and D.~Ruppert, ``Linear non-gaussian component
  analysis via maximum likelihood,'' \emph{arXiv:1511.01609v2}, 2016.

\bibitem{Virta:2016}
J.~Virta, K.~Nordhausen, and H.~Oja, ``Projection pursuit for non-gaussian
  independent components,'' \emph{arXiv preprint arXiv:1612.05445}, 2016.

\bibitem{Theis:2011}
F.~J. Theis, M.~Kawanabe, and K.~R. M\"uller, ``Uniqueness of
  non-gaussianity-based dimension reduction,'' \emph{IEEE Transactions on
  Signal Processing}, vol.~59, no.~9, pp. 4478--4482, 2011.

\bibitem{Blanchard05}
G.~Blanchard, M.~Sugiyama, M.~Kawanabe, V.~Spokoiny, and K.-R. M{\"u}ller,
  ``Non-gaussian component analysis: a semi-parametric framework for linear
  dimension reduction,'' in \emph{Advances in Neural Information Processing
  Systems}, 2005, pp. 131--138.

\bibitem{Bean2014}
D.~M. Bean, ``Non-gaussian component analysis,'' Ph.D. dissertation, University
  of California, Berkeley, 2014.

\bibitem{SasakiNiuSugiyama2016}
H.~Sasaki, G.~Niu, and M.~Sugiyama, ``Non-gaussian component analysis with
  log-density gradient estimation,'' in \emph{Proceedings of the 19th
  International Conference on Artificial Intelligence and Statistics}, 2016,
  pp. 1177--1185.

\bibitem{Cardoso89}
J.-F. Cardoso, ``Source separation using higher order moments,'' in
  \emph{International Conference on Acoustics, Speech, and Signal Processing,
  1989. ICASSP-89.}\hskip 1em plus 0.5em minus 0.4em\relax IEEE, 1989, pp.
  2109--2112.

\bibitem{Nordhausen:2016}
K.~Nordhausen, H.~Oja, and D.~Tyler, ``Asymptotic and bootstrap tests for
  subspace dimension,'' \emph{arXiv:1611.04908}, 2016.

\bibitem{Kankainen07}
A.~Kankainen, S.~Taskinen, and H.~Oja, ``Tests of multinormality based on
  location vectors and scatter matrices,'' \emph{Statistical Methods and
  Applications}, vol.~16, no.~3, pp. 357--379, 2007.

\bibitem{R}
\BIBentryALTinterwordspacing
{R Core Team}, \emph{R: A Language and Environment for Statistical Computing},
  R Foundation for Statistical Computing, Vienna, Austria, 2015. [Online].
  Available: \url{https://www.R-project.org/}
\BIBentrySTDinterwordspacing

\bibitem{ICtest}
\BIBentryALTinterwordspacing
K.~Nordhausen, H.~Oja, D.~E. Tyler, and J.~Virta, \emph{{ICtest}: Estimating
  and Testing the Number of Interesting Components in Linear Dimension
  Reduction}, 2016, {R} package version 0.2. [Online]. Available:
  \url{https://CRAN.R-project.org/package=ICtest}
\BIBentrySTDinterwordspacing

\bibitem{YeWeiss:2003}
Z.~Ye and B.~Weiss, ``Using the bootstrap to select one of a new class of
  dimension reduction methods,'' \emph{Journal of the American Statistical
  Association}, vol.~98, no. 464, pp. 968--979, 2003.

\bibitem{LuoLi:2016}
W.~Luo and B.~Li, ``Combining eigenvalues and variation of eigenvectors for
  order determination,'' \emph{Biometrika}, vol. 103, no.~4, pp. 875--887,
  2016.

\bibitem{Oja:2006}
H.~Oja, S.~Sirki\"a, and J.~Eriksson, ``Scatter matrices and independent
  component analysis,'' \emph{Austrian Journal of Statistics}, vol.~35, pp.
  175--189, 2006.

\bibitem{TCDO09}
D.~Tyler, F.~Critchley, L.~D\"umbgen, and H.~Oja, ``Invariant coordinate
  selection,'' \emph{Journal of Royal Statistical Society, Series B}, vol.~71,
  pp. 549--592, 2009.

\bibitem{EatonTyler91}
M.~L. Eaton and D.~E. Tyler, ``On {W}ielandt's inequality and its application
  to the asymptotic distribution of the eigenvalues of a random symmetric
  matrix,'' \emph{The Annals of Statistics}, vol.~19, no.~1, pp. 260--271,
  1991.

\end{thebibliography}
%



%







\end{document}